\documentclass[a4paper,12pt, reqno]{amsart}

\usepackage{fullpage}
\usepackage{graphicx}
\usepackage{hyperref}
\usepackage{amsmath}
\usepackage{amssymb}
\usepackage{amsthm}
\usepackage{appendix}
\usepackage{verbatim}

\usepackage[left=2.5cm,right=2.5cm,top=2.5cm,bottom=2.5cm,a4paper,footskip=1cm]{geometry}

\makeatletter
\renewcommand{\pod}[1]{\mathchoice
  {\allowbreak \if@display \mkern 18mu\else \mkern 8mu\fi (#1)}
  {\allowbreak \if@display \mkern 18mu\else \mkern 8mu\fi (#1)}
  {\mkern4mu(#1)}
  {\mkern4mu(#1)}
}

\makeatletter
\def\old@comma{,}
\catcode`\,=13
\def,{%
  \ifmmode%
    \old@comma\discretionary{}{}{}%
  \else%
    \old@comma%
  \fi%
}
\makeatother

\allowdisplaybreaks[1]

\theoremstyle{plain}
\newtheorem{thm}{Theorem}
\newtheorem{lem}{Lemma}

\newtheorem{conj}{Conjecture}
\newtheorem{cor}{Corollary}
\theoremstyle{definition}

\begin{document}
\title{On Polignac's conjecture and arithmetic progressions}
\author{Stijn S.C. Hanson}
\maketitle

\begin{abstract}
  In this paper we investigate the recent advances by Zhang, Maynard and Pintz towards Polignac's conjecture and give some new results concerning the relationship between Polignac numbers and arithmetic progressions.
\end{abstract}

\section{Introduction}
We start with a statement of Polignac's conjecture:
\begin{conj}[de Polignac, 1849]\label{polignac}
  Let $k$ be any positive even integer and let $p_n$ be the n\textsuperscript{th} prime number. Then, for infinitely many $n \in \mathbb{N}$, we have $p_{n+1} - p_n = k$.\cite{polignac}
\end{conj}
Any $k$ that satisfies Polignac's conjecture is called a Polignac number and the twin prime conjecture is clearly equivalent to asking whether 2 is a Polignac number.

This conjecture was partially resolved in \cite{maynard2} and \cite{zhang} with the following theorem
\begin{thm}[Maynard-Zhang, 2013]\label{mayzha}
  There exist infinitely many pairs of prime numbers whose difference is no more than $600$.
\end{thm}
This clearly shows that there exists at least one Polignac number $2 \leq k \leq 600$ and this, together with the following theorem from J\'{a}nos Pintz\cite{pintz}, asserts that there are infinitely many Polignac numbers.
\begin{thm}[Pintz, 2013]\label{pintz}
  There is an ineffective constant $C$ such that every interval of the form $[m, m + C]$ contains a Polignac number.
\end{thm}
Unfortunately this theorem is ineffective - it offers no way to generate future Polignac numbers - but it does give direction for future research.

\section{Results}
Let $\mathcal{H} = \{h_1, h_2, \dots, h_k\}$ be a set of integers. $\mathcal{H}$ is said to be admissible if, for every prime number $p$, there is some $n \in \mathbb{N}$ such that, for every $1 \leq i \leq k$
\begin{equation}
  h_i \not\equiv n \pmod{p}
\end{equation}

The way that Theorem \ref{mayzha} was proven was showing that, if $|\mathcal{H}| \geq k$ is admissible, then $n + \mathcal{H}$ contains two primes infinitely often. Our current lowest value for $k$ is $59$ \cite{polymathbounds}.

If we suppose that we have proved that any admissible set of size $k$ has infinitely many translates that contain $\rho$ primes then we can verify Conjecture \ref{polignac} in a few specific cases. 

\begin{samepage}
  \begin{lem}\label{han1}
    Let $d = \prod_{p \leq k} p$. Then, for every $N \in \mathbb{N}$,
    \begin{equation}\label{adm}
      \{dN, 2dN, \dots, (k - 1)dN\} 
    \end{equation}
    contains at least one Polignac number.
    \begin{proof}
      Consider the set
      \begin{equation}
        \mathcal{H} = \{0, dN, 2dN, \dots, (k-1)dN\}.
      \end{equation}
      If we take any prime $p \leq k$ then, for all $h \in \mathcal{H}$, we have $p | h$. I.E.
      \begin{equation}
        h \equiv 0 \pmod{p}
      \end{equation}
      for all $p \leq k$. Furthermore, if $p > k$ then $\mathcal{H}$ cannot occupy all residue classes modulo $p$ as there are only $k$ elements. Therefore $\mathcal{H}$ is admissible. 
      
      As $k$ is large enough we know that there are infinitely many translates of $\mathcal{H}$ which contain two primes infinitely often. Therefore there is some pair $(n_1 dN, n_2 dN)$ which are simultaneously prime infinitely often, and the difference between these must be one of the elements of the set (\ref{adm}).
    \end{proof}
  \end{lem}
\end{samepage}
Let $q \in \mathbb{N}$ and consider the sequence given by $n \equiv 0 \pmod{q}$. This has a subsequence given by $n \equiv 0 \pmod{qd}$ where $d = k\#$ is the primorial of $k$. This, in turn, has the subsequence
\begin{equation}
  (qd, 2qd, \dots, (k-1)qd), (kqd, 2kqd, \dots, (k-1)kqd), \dots
\end{equation}
where each of the parenthesised terms contains a Polignac number by the above lemma so we have proven:
\begin{cor}
  There are infinitely many Polignac numbers on any arithmetic progression of the form
  \begin{equation}
    q, 2q, \dots
  \end{equation}
  for all $q \in \mathbb{N}$.
\end{cor}

We can extend this result fairly easily to give a Dirichlet-type theorem.
\begin{thm}
  Every arithmetic progression $a, a + q, a + 2q, \dots$ such that $q \mid a$ contains infinitely many Polignac numbers.
  \begin{proof}
    We have just proved the case where $a = 0$ so suppose that $a \neq 0$. Now define
    \begin{equation}
      N_i = \frac{a}{q}\left(i \prod_{p \leq k } p - 1\right).
    \end{equation}
    Then
    \begin{equation}
      a + N_iq = a i \prod_{p \leq k } p
    \end{equation}
    is a subsequence of our original arithmetic progression. But every parenthesised term
    \begin{equation}
      (a + N_1 q, a + N_2 q, \dots a + N_{k-1}q) + (a + N_k q, a + 2 N_k q , \dots, a + N_{k-1}N_k q), \dots
    \end{equation}
    contains a Polignac number and so the whole arithmetic progression must contain infinitely many.
  \end{proof}
\end{thm}

This raises the question as to whether the following statement is true.
\begin{conj}
  Every non-trivial arithmetic progression contains infinitely many Polignac numbers where the trivial arithmetic progressions are those which are entirely odd.
\end{conj}

Lemma \ref{han1} implies the infinitude of Polignac numbers but in a more effective form than Theorem \ref{pintz}.

By using Pintz' result that every interval of the form $[m, m + C]$ contains a Polignac number if $C$ is large enough we can calculate a lower bound for the upper asymptotic density of the set of Polignac numbers $\mathcal{P}$.
\begin{align}
  \overline{\sigma}(\mathcal{P}) &= \limsup_{n \rightarrow \infty} \frac{|\mathcal{P} \cap [0, n]|}{n} \nonumber \\
  &\geq \limsup_{n \rightarrow \infty}\frac{\left|\bigcup_{i=1} ^{n/c}\mathcal{P}\cap [(i-1)C, iC]\right|}{n}.
\end{align}
All of the intersections above are non-empty but they might count some Polignac numbers twice so this is
\begin{align}
  &\geq \limsup_{n \rightarrow \infty}\frac{\sum_{i=1} ^{n/c} 1}{2n} \nonumber \\
  &\geq \limsup_{n \rightarrow \infty} \frac{n/c - 1}{2n} \nonumber \\
  &= \frac{1}{2c}.
\end{align}

This is positive and so we can apply Szemer\'{e}di's theorem \cite{szemeredi} \cite{taovu} to show
\begin{thm}
  The set of Polignac number contains arbitrarily long arithmetic progressions.
\end{thm}

Which has the following natural generalisation.

\begin{conj}
  The Polignac numbers contain infinitely long arithmetic progressions.
\end{conj}
\section*{Acknowledgements}
The author would like to thank Martina Balagovic, Samuel Bourne, Alan Haynes, James Maynard, Terence Tao, Robert Vaughan and Charles Wing for many useful conversations and corrections. Further thanks must go to Christopher Hughes for supervising the project which this was borne out of.

\end{document}